\documentclass[12pt]{article}
\usepackage{latexsym}
\usepackage{amssymb}
\usepackage[dvips]{graphicx}
\usepackage{epsfig}
\usepackage{rotating}
\usepackage{amsfonts}
\usepackage{amsmath}

\usepackage{hyperref}
\usepackage{pstricks}

\newtheorem{lemma}{Lemma}[section]
\newtheorem{theorem}{Theorem}[section]
\newtheorem{corollary}{Corollary}[section]
\newtheorem{definition}{Definition}[section]
\newtheorem{proposition}{Proposition}[section]

\newcommand{\beqa}{\begin{eqnarray}}
\newcommand{\eeqa}{\end{eqnarray}}

\newcommand{\bea}{\begin{eqnarray}}
\newcommand{\eea}{\end{eqnarray}}

\newcommand{\be}{\begin{equation}}
\newcommand{\ee}{\end{equation}}
\newcommand{\beq}{\begin{eqnarray}}
\newcommand{\eeq}{\end{eqnarray}}

\def\lsim{\
  \lower-2.0pt\vbox{\hbox{\rlap{$<$}\lower5.5pt\vbox{\hbox{$\sim$}}}}\ }
\def\gsim{\
  \lower-2.0pt\vbox{\hbox{\rlap{$>$}\lower5.5pt\vbox{\hbox{$\sim$}}}}\ }

\begin{document}

\begin{titlepage}
\begin{flushright}
\end{flushright}

\vspace{20pt}

\begin{center}

{\Large\bf 
Multi-orientable Group Field Theory
}
\vspace{20pt}

 Adrian Tanasa

\vspace{15pt}

\vspace{10pt}$^{}${\sl
LIPN, Institut Galil\'ee, CNRS UMR 7030,\\
Univ. Paris Nord, 99 av. Cl\'ement, 93430 Villetaneuse, France, UE}\\

\vspace{10pt}
$^{}${\sl
Horia Hulubei National Institute for Physics and Nuclear Engineering,\\
P.O.B. MG-6, 077125 Magurele, Romania, EU
}\\

\vspace{20pt}
E-mail:  
 {\em adrian.tanasa@ens-lyon.org}

\vspace{10pt}

\begin{abstract}
\noindent
Group Field Theories (GFT) are quantum field theories over group manifolds; they can be seen as a generalization of matrix models.
GFT Feynman graphs are tensor graphs 
generalizing ribbon graphs (or combinatorial maps); these graphs are 
dual not only to manifolds. In order to simplify the topological structure of 
these various singularities, colored GFT was recently introduced and intensively studied since. We propose here a different simplification of GFT, which we call {\it multi-orientable} GFT.

We study the relation between multi-orientable GFT Feynman graphs and colorable graphs. We prove that tadfaces and some generalized tadpoles are absent. 
Some Feynman amplitude computations 
are performed. A few remarks on the renormalizability of both multi-orientable and colorable GFT are made. A generalization from three-dimensional to four-dimensional theories is also proposed.
\end{abstract}

\end{center}


\end{titlepage}


\setcounter{footnote}{0}

\section{Introduction and motivation}
\label{Intro}
\renewcommand{\theequation}{\thesection.\arabic{equation}}
\setcounter{equation}{0}

Graph theory \cite{graph-books} is known to play a fundamental r\^ole in describing the combinatorics of quantum field theory (QFT) \cite{qft2, books-QFT, io-SLC}. A natural generalization of graphs, the combinatorial maps \cite{maps} (or ribbon graphs) play the same r\^ole for matrix models, known to be related to non-commutative QFT \cite{ncqft} or to two-dimensional quantum gravity \cite{2d, 2d2}.

The generalization of these models to three- and four-dimensional space is known under the name of group field theory (GFT) \cite{gft}. Feynman amplitudes of GFT can be equivalently expressed as spin-foam models \cite{sf}, which together with string theory, represent one of the most developed approaches today for a fundamental theory of quantum gravity \cite{dia}.

Nowadays, GFT is allocated a great deal of interest by the mathematical physics community
(see  \cite{gft-regain} and references within). 
The r\^ole that combinatorial maps play for matrix models is now played by tensor graphs (these graphs have three strands per edge in three-dimensional models and respectively four strands per edge in four-dimensional models).

Nevertheless, these tensor graphs are dual not only to manifolds; 
the topological structure of their singularities is indeed very complicated.
This can be seen as a drawback or, on the contrary, as a mathematical richness (for example, in QFT models on the non-commutative Moyal space, it is the ``non-planar'' tadpole-like graphs - definitely more complicated from a topological point of view that the planar ones - that give rise to the celebrated phenomenon of ultraviolet/infrared mixing \cite{uv/ir}; getting rid of these ``non-planar'' graphs would eliminate this new, physically rich phenomenon). 

For the sake of completeness, let us mention here that other ways of implementing noncommutative QFT exist in the literature, ways which do not lead to the ultraviolet/infrared mixing (see for example \cite{bahns}). 
Taking yet another point of view, one can use braided QFT and this also leads to no ultraviolet/infrared mixing \cite{braided}. 


\medskip

Coming back to GFT, in order to eliminate a whole class of ``wrapping singularities'', a restricted class of models was introduced - the colored GFT
\cite{color}. The idea is to have each of the $D$ edges incoming/outgoing to a vertex colored with a distinct color ($D$ being the dimension of space-time) and to have two types of vertices. From a graph theoretical point of view, this last part 
means that one restricts the set of graphs to {\it bipartite graphs} (also known as {\it bigraphs} - see for example \cite{bipartite}). Related to this topic of colored GFT, let us also mention here the recent literature on manifold crystallization \cite{cristal}. Moreover, it is this notion of bigraphs which plays a key-r\^ole when proving that the color in colored GFT guarantees orientability of the piecewise linear pseudo-manifold associated to each graph of the perturbative expansion \cite{cara}.

Colored GFT has been  given lately a great deal of interest 
within GFT;  
several 
achievement concerning 
these colored models have been made (see \cite{dev-cgft}, \cite{bgo}, \cite{ryan}, \cite{j2}, \cite{tf} and references within).

\medskip

In this paper, we propose a simplification of the GFT tensor graph class which is different from the coloring one; we call this new simplification {\it multi-orientability}. The idea behind this is to introduce a notion of orientability at the level of the GFT vertex. This was already successfully implemented for ribbon graphs in QFT on Moyal space; the respective non-commutative models were intensively studied in various papers of mathematical physics (see \cite{o1}, \cite{orient-ncqft} and references within).

We introduce this new class of GFT models both in three and in four dimensions. We also analyze the relation between the multi-orientable GFT graphs and the colorable ones. 
We prove that tadfaces and some class of generalized tadpoles are not allowed within this multi-orientable framework. Our proofs also allows to recover the known result that tadfaces are also forbidden within the colorable framework.
Furthermore, a particular class of multi-orientable GFT graphs is identified. The associated Feynman amplitudes for the BF-theories are computed, both in three (the Boulatov model \cite{boulatov}) and in four dimensions (the Ooguri model \cite{ooguri}).

Moreover, we illustrate the following phenomenon. Within the framework of colorable GFT, divergent two- and four-point graphs (for the three-dimensional case) represent quantum corrections of types not present in the bare action of the model. Within the multi-orientable framework, these graphs represent quantum corrections of type already existing in the bare action.

\medskip

The paper is organized as follows. In the following section we present the definition of colorable Boulatov three-dimensional model and we give the definition of multi-orientable GFT. The following section deals with graph theoretical consequences of this definition. We look closely to the issues of tadpoles and tadfaces and we establish the relation between colorable and multi-orientable three-dimensional GFT  graphs. In the fourth section we compute some Feynman amplitudes of non-colorable but multi-orientable, colorable (and multi-orientable) and also non-multi-orientable (and non-colorable) graphs. The following section presents some considerations with respect to the renormalizability of the colorable and of the multi-orientable models. The sixth section proposes a generalization of multi-orientability to the four-dimensional case. The last section presents some perspectives for future developments of this new type of models.

\section{Multi-orientable GFT in three dimensions}
\label{sec:3d}
\renewcommand{\theequation}{\thesection.\arabic{equation}}
\setcounter{equation}{0}

In this section we first recall the definition of the colorable models and we then give the definition of the multi-orientable ones.

The GFT models (colorable or non-colorable) usually studied in literature are {\it orientable}: this means that the simplicial complexes dual to these graphs are orientable. 
We keep this simplification in this paper.

The field $\phi$ of the three-dimensional topological Boulatov model is a map
$\phi:G^3\to \mathbb R$, 
where the group $G$ is the $SU(2)$ one.
For such an orientable model, the action  writes
\beqa
\label{act-boul0}
&&S[\phi]=\frac 12 \int dg_1 dg_2 dg_3 \phi(g_1,g_2,g_3)\phi(g_1,g_2,g_3)\\
&&+\frac{\lambda}{4!}
\int dg_1\ldots dg_6 dg_{1'}\ldots dg_{6'} \phi (g_1,g_2,g_3)\phi(g_{3'},g_4,g_5)\phi(g_{5'},g_{2'},g_6)\phi(g_{6'},g_{4'},g_{1'})\nonumber\\
&&\delta(g_1g_{1'}^{-1})\ldots \delta(g_6g_{6'}^{-1}),\nonumber
\eeqa
which is equivalent to the simplified form:
\beqa
\label{act-boul}
S[\phi]&=&\frac 12 \int dg_1 dg_2 dg_3 \phi(g_1,g_2,g_3)\phi(g_1,g_2,g_3)\\
&+&\frac{\lambda}{4!}\int dg_1\ldots dg_6 \phi (g_1,g_2,g_3)\phi(g_3,g_4,g_5)\phi(g_5,g_2,g_6)\phi(g_6,g_4,g_1),\nonumber
\eeqa
Note that the integration over the group is done with the invariant Haar measure, both in \eqref{act-boul0} and in \eqref{act-boul}.

Let us make here the following remark. The integral kernel represented by the product of $\delta$ functions in \eqref{act-boul0} can be interpreted as a natural generalization, within the GFT formalism, of the crucial notion of {\it locality} of QFT. Thus, if one takes the example of the local $\Phi^4$ model, the non-quadratic part (the interaction) of the action in configuration (or direct) space writes:
\beqa
\label{local}
 \int d x d x_1\ldots d x_4 \delta (x-x_1)\ldots \delta(x-x_4)\Phi(x_1)\ldots \Phi(x_4)=\int d x \Phi^4(x).
\eeqa
It is this form of the vertex which is generalized in formula \eqref{act-boul0} and \eqref{act-boul}. Later on, within the framework of QFT on the non-commutative Moyal space, this notion of locality was replaced by a more adapted notion of ``Moyality'' (see again \cite{ncqft} and references within). We thus advocate to keep the general form \eqref{act-boul0} for the vertex of {\it any} GFT model, topological or not, in three or in four dimensions, in order to have a natural generalization of the notion of locality (and of ``Moyality''). Let us stress that this has been also defended for example in \cite{io-5}.

Nevertheless, a  different point of view is usually taken in the GFT literature, point of view inspired from the spin-foam experience in quantum gravity. Thus, for a non-topological model, the interaction is the one which is modified (and not the propagator), because in this way one has a clear geometrical meaning of the added constraints, already at the level of the action. However, in order to perform 
renormalizability studies of GFT, it seems natural to us to switch to the QFT point of view advocated here, for the reasons explained above. This idea can be defended also by the fact that, in theoretical physics, renormalization is done within some QFT framework.

\medskip

\begin{figure}
\begin{center}
\includegraphics[scale=0.2,angle=0]{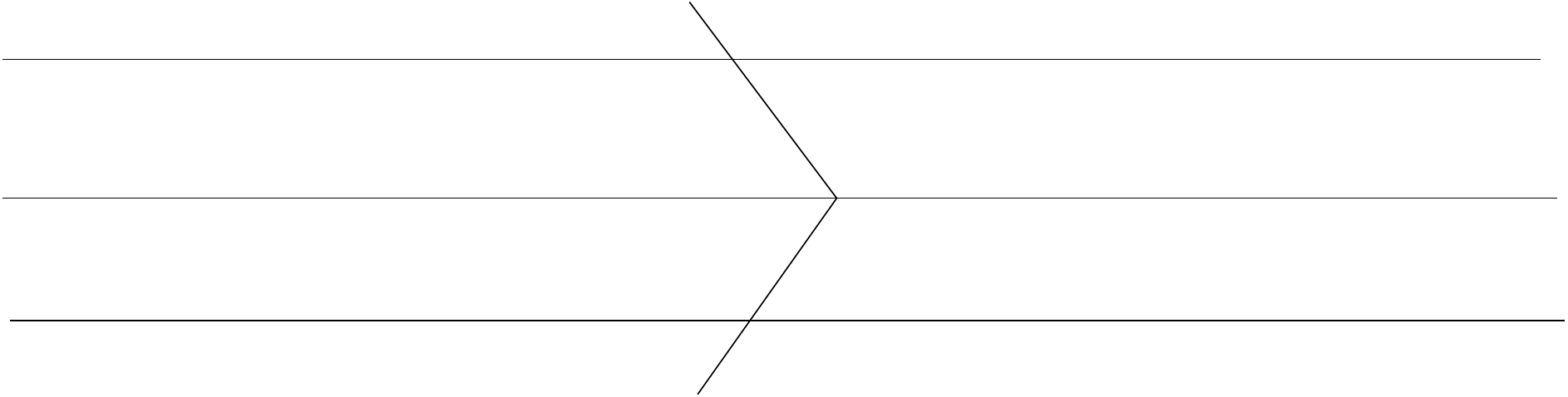}
\caption{The edge of the orientable GFT tensor graphs.}
\label{fig:edge}
\end{center}
\end{figure}

\begin{figure}
\begin{center}
\includegraphics[scale=0.3,angle=0]{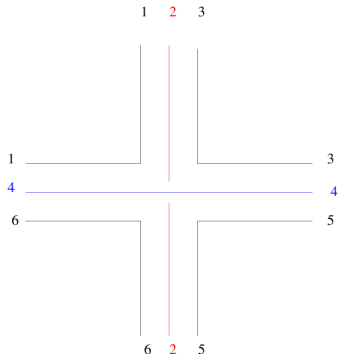}
\caption{The vertex of the Boulatov model considered here.}
\label{fig:vertex}
\end{center}
\end{figure}

The Boulatov field $\phi$, taken real-valued, is not assumed here to have specific symmetry properties under the permutations of its arguments. Note that
the more involved case, where not only the identical permutation is associated to such an edge, is analyzed in \cite{not-ori} or in \cite{not-ori2}.
Thus, the edge of the graphs considered here is represented like in Fig. \ref{fig:edge}.
The propagator of the model writes
\beqa
\int dh \prod_{i=1}^3 \delta(g_ihg_{i'}^{-1}).
\eeqa
The form of the vertex chosen in the action \eqref{act-boul} corresponds to the one of Fig. \ref{fig:vertex}.

The colorable three-dimensional model is defined in the following way. The (real-valued) Boulatov field is replaced by four complex-valued field $\phi_p$, the index $p=0,\ldots,3$ being referred to as some color index. Moreover, one now has two types of interactions, a $\phi^4$ one and a $\bar \phi^4$ one. Furthermore, a clockwise cyclic ordering at one of the types of vertices (and an anticlockwise at the second type of vertex) of the four colors at the vertex is imposed; thus the action writes
\beqa
\label{act:col}
S_{{\mathrm{col}}}=\frac 12 \sum_p \int \bar \phi_p \phi_p + \frac{\lambda}{4!}\int \phi_0 \ldots \phi_4 +
\frac{\lambda}{4!}\int \bar \phi_0 \ldots \bar \phi_4, 
\eeqa
where the integrations over the group are left implicit.

From a graph-theoretical point of view, this means that one imposes a four-coloring of the edges and that only bipartite graphs are kept. Furthermore, one has to respect the cyclic ordering of the two types of vertices described above.

\medskip

\begin{figure}
\begin{center}
\includegraphics[scale=0.15,angle=0]{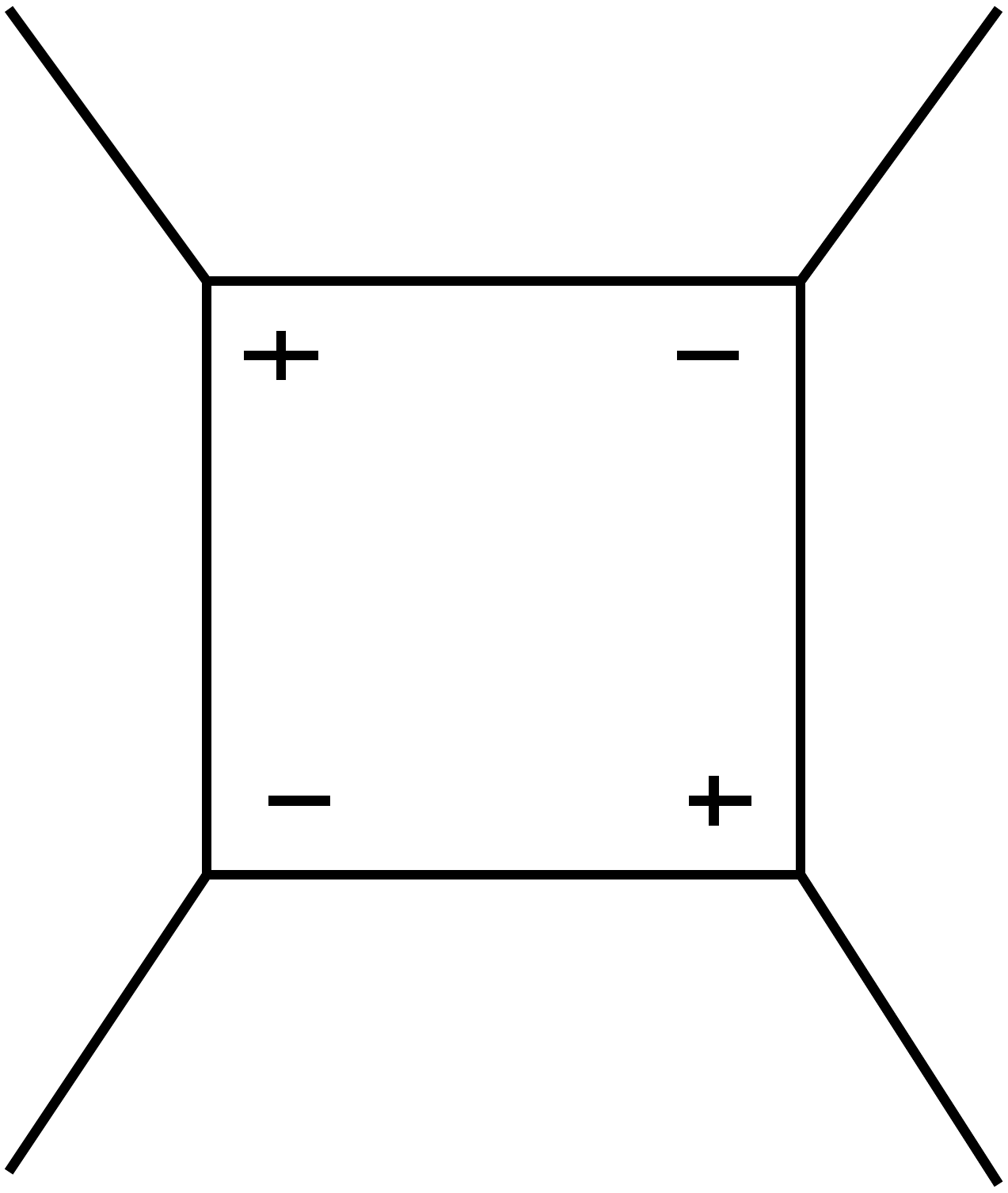}
\caption{The vertex of the proposed multi-orientable GFT model.}
\label{vertex-o}
\end{center}
\end{figure}

Let us now introduce the announced three-dimensional multi-orientable model. 
As in the colored case, one has a complex-valued field $\phi$. Nevertheless, we do not need more copies of this field.
The Boulatov interaction is restricted to vertices where each corner has a $+$ or a $-$ label. Furthermore, each vertex has two corners labeled with $+$ and two corners labeled with $-$, which are cyclically ordered as shown in Fig. \ref{vertex-o}. 
A field propagates from a $-$ to a $+$ corner of some vertex.

This notion of corners of a vertex is natural within the framework of non-local QFT (as is the Moyal non-commutative QFT, where, as mentioned in the Introduction, the idea of multi-orientability has proved to be very useful). Nevertheless, GFT also can be seen in some sense as a non-local QFT on the group manifold on which it lives, since the interaction is not defined on some ``group point'' $g$ (see \eqref{act-boul0} and \eqref{act-boul}), unlike the interaction of the local QFT models (see \eqref{local}).
  
We propose to call this model {\it multi-orientable}, because one has:
\begin{itemize}
\item on the one hand, the usual orientability of the GFT propagation (see Fig. \ref{fig:edge} above)
\item and on the other hand, the orientability of the vertex.
\end{itemize}

The action of the model writes
\beqa
\label{act:mo}
S[\phi]=\frac 12 \int \bar \phi \phi + \frac{\lambda}{4!}\int \bar \phi \phi \bar \phi \phi,
\eeqa
where, as in \eqref{act:col}, the integrations over the group are left implicit.

At a graph-theoretical level, the model \eqref{act:mo} generates Feynman graphs built up of vertices like the one represented in Fig. \ref{vertex-o} and edges like the one 
represented in Fig. \ref{fig:edge}.

\medskip

For the sake of completeness, let us also mention that combinatorial maps with four-valent vertices like the orientable ones of Fig. \ref{vertex-o} have not been counted \cite{fusy}; this is not the case for maps (of genus $g$) in general, which are very well analyzed from a combinatorial point of view (see again \cite{maps} and references within).

\medskip

Before ending this section, let us also remark that some kind of hybrid model, where the vertex is kept orientable 
but in which one not only
allows the identical permutation of the three arguments of the edge, can be naturally defined.

\section{Multi-orientable graphs; relation with the colorable ones}
\renewcommand{\theequation}{\thesection.\arabic{equation}}
\setcounter{equation}{0}

\begin{figure}
\begin{center}
\includegraphics[scale=0.15,angle=0]{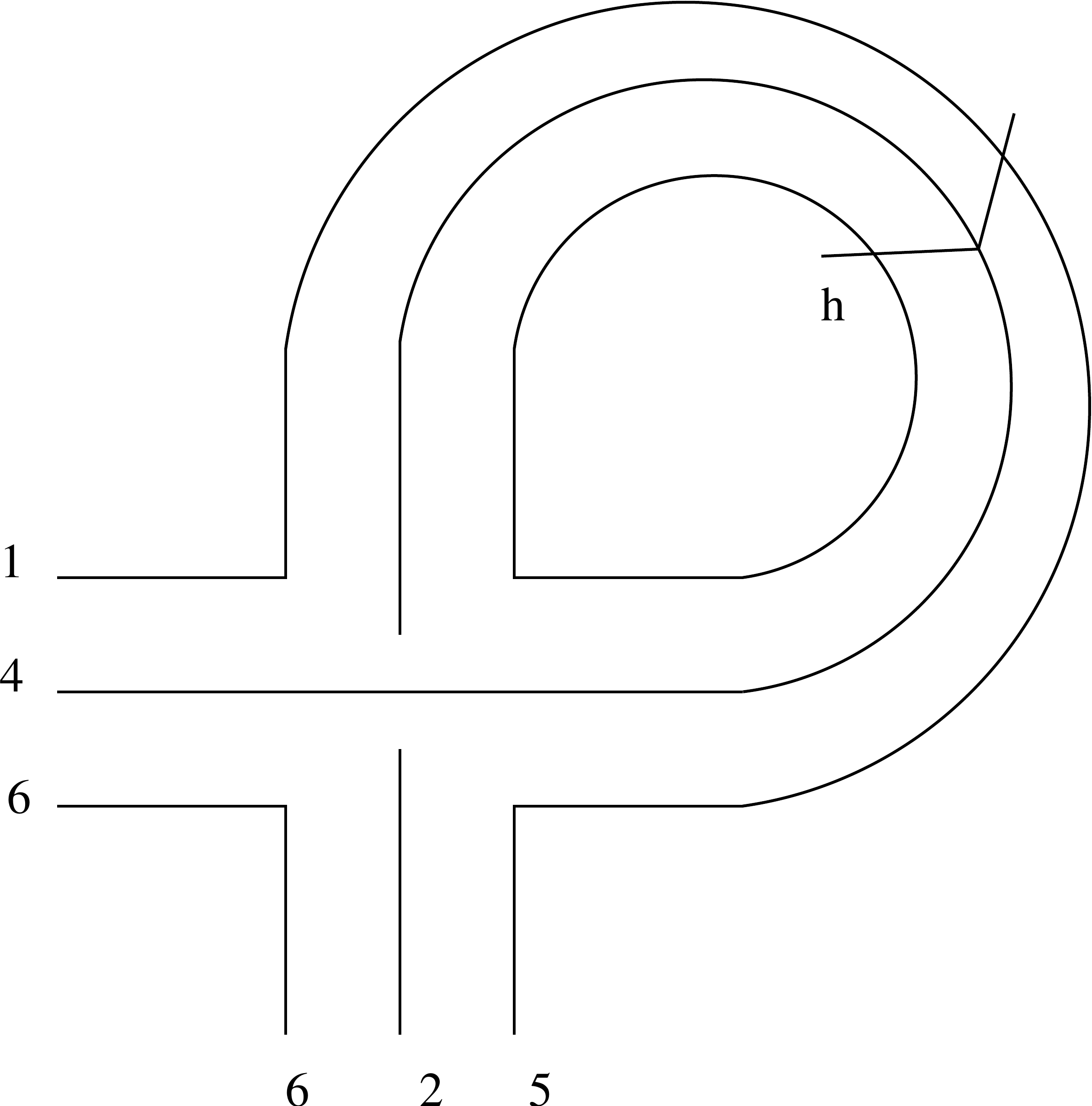}
\caption{A tadpole graph which is multi-orientable. This graph is not allowed by the colorable GFT. The group elements of the external edges are $g_1$, $g_4$, $g_6$ , $g_2$ and $g_5$ while the group element $h$ is associated to the single internal edge.}
\label{tad1}
\end{center}
\end{figure}

\begin{figure}
\begin{center}
\includegraphics[scale=0.15,angle=0]{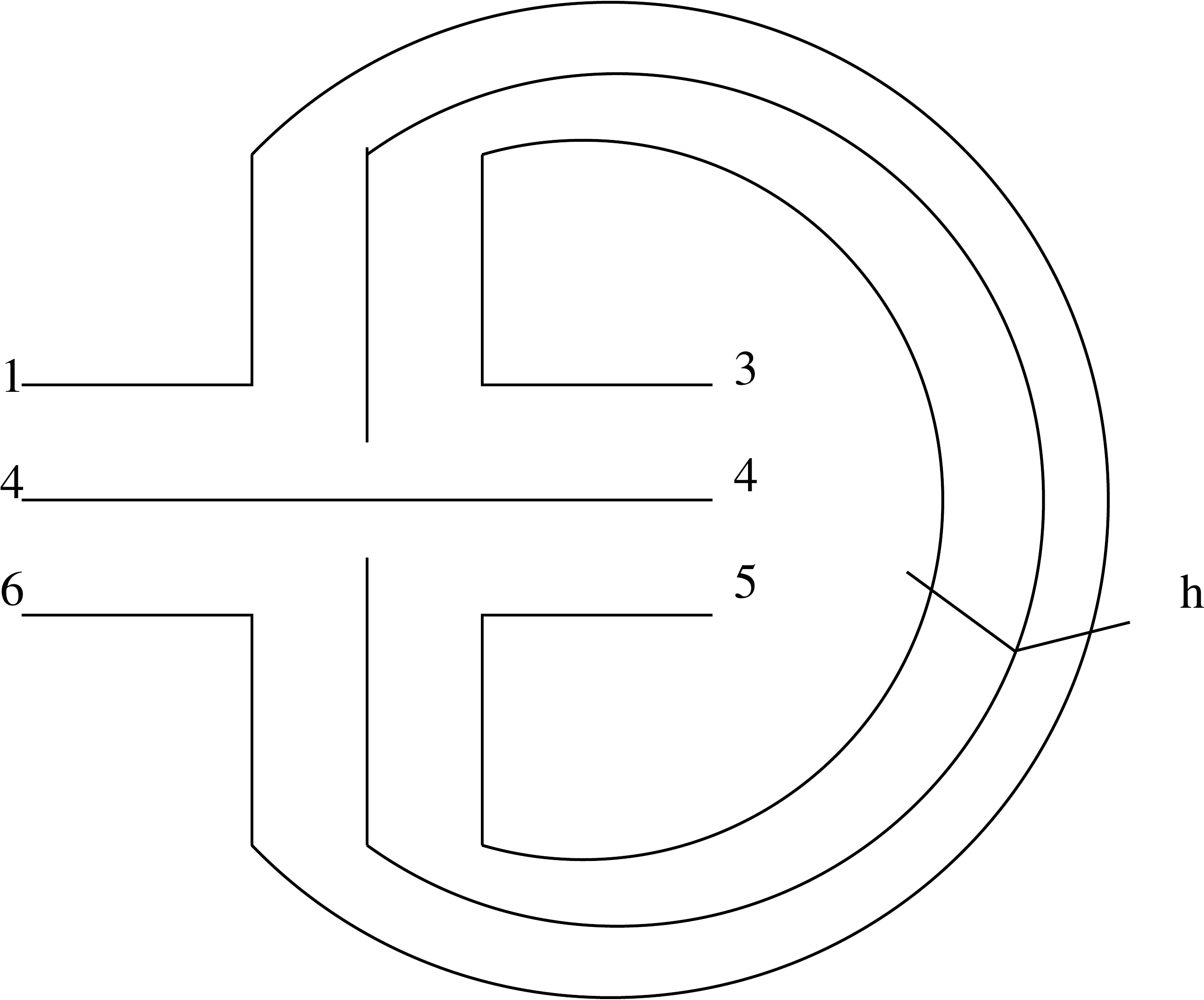}
\caption{A tadpole graph which is not multi-orientable. This type of tadpole is thus discarded within the multi-orientable GFT framework. The group elements of the external edges are $g_1$, $g_4$, $g_6$ , $g_3$ and $g_5$ while the group element $h$ is associated to the single internal edge.}
\label{tad2}
\end{center}
\end{figure}

Let us take a closer look at the graph-theoretical consequences of 
the definition of the multi-orientable models of the previous section.
The colorability discards a highly significant class of graphs, including the so-called ``wrapping singularities'', which correspond to graphs containing tadpoles like Fig. \ref{tad1} 
and \ref{tad2} 
(see again \cite{color} for more details).


Let us now prove the following result:

\begin{proposition}
\label{inclu}
Every GFT graph which is colorable is also multi-orientable.
\end{proposition}
{\it Proof.} Recall that the conditions imposed for a graph to be orientable are: the existence of two types of vertices, cyclic ordering (clockwise and anti-clockwise) at these vertices and the fact that the four edges adjacent to such a label wear a specific label (the color). The conditions to define multi-orientability are actually just a subset of these conditions. Thus, the cyclic ordering is imposed, but no labeling (coloring) of the edges - there is no distinction between the two $+$ or the two $-$ corners of a vertex. Moreover, the proof is completed by the following fact. If the multi-orientability condition is satisfied as some vertex, then it is also satisfied at another vertex connected to the first one by an edge of some generic color. (QED)

\medskip

\begin{figure}
\begin{center}
\includegraphics[scale=0.15,angle=0]{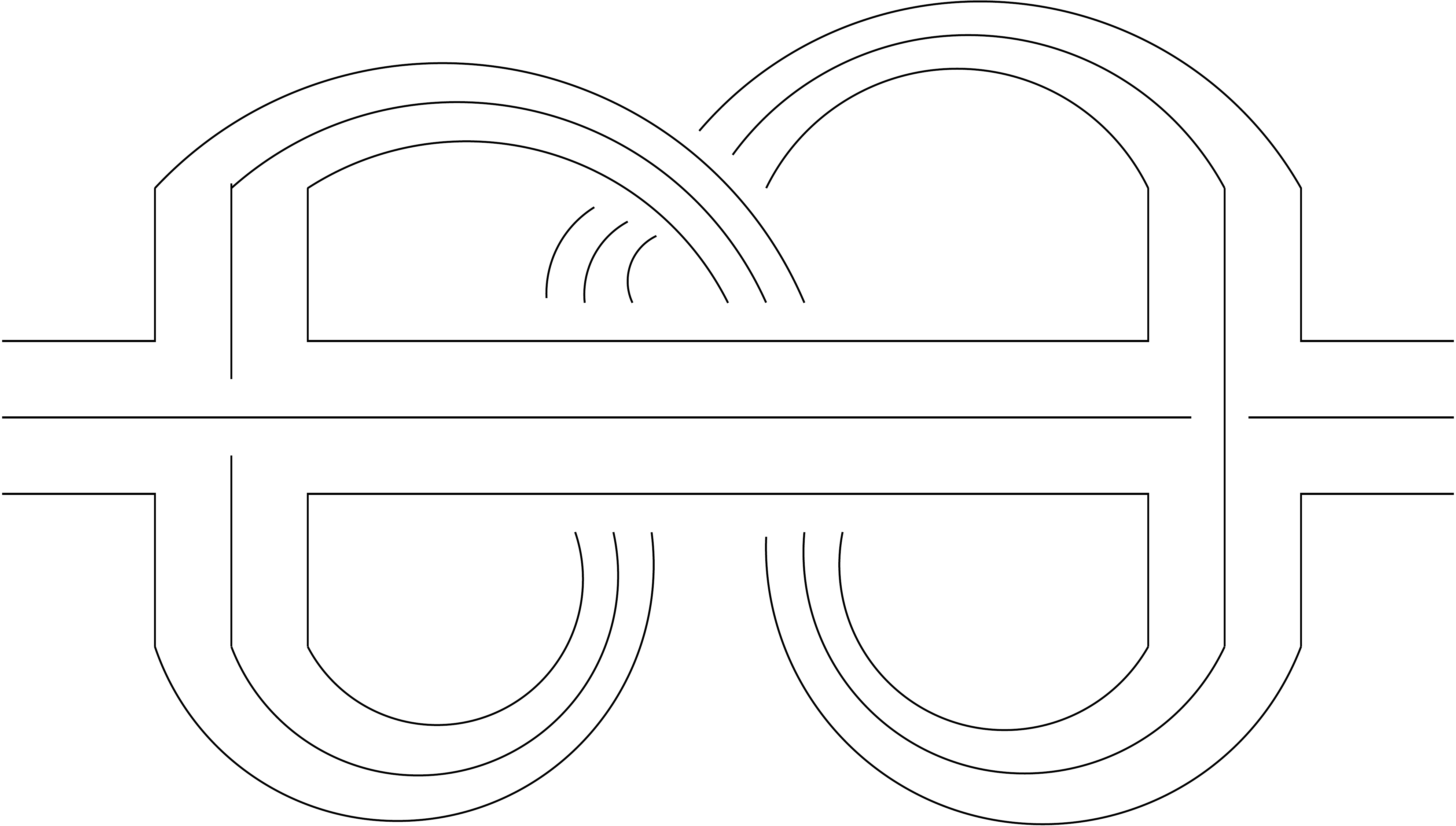}
\caption{A non-tadpole example of a GFT graph which is non-colorable but is multi-orientable.}
\label{tadface}
\end{center}
\end{figure}

The reciprocal statement is not true. A counterexample is the tadpole graph of Fig. \ref{tad1}, which is multi-orientable (or, in other words, it can {\it multi-oriented}) but is non-colorable. Another example, which is not a tadpole graph, is the two-point graph of Fig. \ref{tadface}.

\medskip

All of this can be rephrased as: 
multi-orientability discards a less important class of graphs then colorability.

\subsection{Tadpoles and generalized tadpoles}

Before investigating the issue of tadpoles within the multi-orientable framework, let us make the following remark. Since the edges of the models we deal with here (colorable or not, multi-orientable or not) do not allows twists (as already stated in the previous section), one can drop the ``middle'' strand of any edge and obtain an one-to-one correspondence with some ribbon graph (or combinatorial map). The ribbon graph thus obtained is the {\it jacket} introduced in \cite{tadface} and later generalized in \cite{j2}. 
Furthermore, in \cite{ryan}, it was showed that the jacket graphs represent (Heegaard) splitting surfaces for the triangulation dual to the Feynman graph; this
allow to re-express the Boulatov model as a QFT model on these Riemann surfaces (see again \cite{ryan}). 

One can thus refer to the planarity of the respective tensor GFT graph as to the planarity of the ribbon graph associated in this way.
Moreover, one can count the number of faces broken by the external legs; we denote this number by $B$. If $B>1$ we call the respective graph {\it irregular} (see again \cite{ncqft} for details).

The tadpole in Fig. \ref{tad1} is thus referred to as a {\it planar} tadpole, while the tadpole in Fig. \ref{tad2} is referred to as a {\it ``non-planar''} tadpole (although the terminology ``non-planar'' usually used in the literature is incorrect, because the respective ribbon graph is planar but it just has a number of broken faces superior to one).

We have seen above that planar tadpoles are allowed by multi-orientable models, while the non-planar ones are not.

Let us now recall from \cite{tf} the following definition:

\begin{definition}
A generalized tadpole is a graph with one external vertex.
\end{definition}

Planar generalized tadpoles are allowed by multi-orientability. 
A ``non-planar'' generalized tadpole is a graph with two external edges and with $B=2$. These graphs are not allowed by multi-orientability (this is a result already known from the non-commutative QFT literature, see for example \cite{o1}).

\subsection{Tadfaces}


We recall from \cite{tadface} the following definition:

\begin{definition}
\label{d:tf}
A  { tadface} is a face that goes at least
twice through a line.
\end{definition}

\begin{figure}
\begin{center}
\includegraphics[scale=0.1,angle=0]{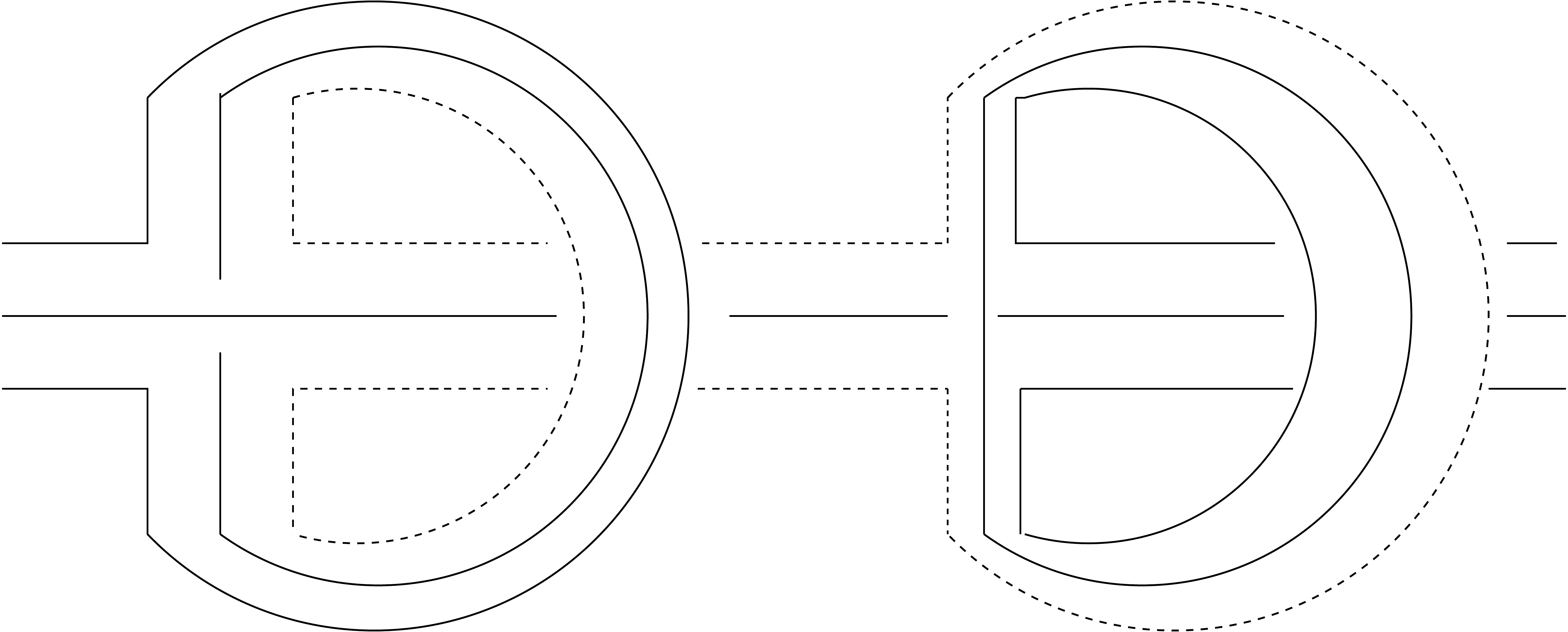}
\caption{A non-multi-orientable tadface graph. The tadface is represented by the dashed line. One notices that the respective face goes twice through the edge relying the two tadpoles, once through one strand and once through another strand of the edge. }
\label{tf1}
\end{center}
\end{figure}

Let us give a few more explanation of this. Such a tadface can be obtained if, 
one goes through the respective edge a first time through a strand and a second time through a second strand. A second-order example of such a graph is given in Fig. \ref{tf1} (graph which is to be inserted into a ``bigger'' graph, such that the usual 1PI condition can be kept). Nevertheless, this graph is not multi-orientable, since it is made up of two non-multi-orientable tadpoles like the ones in Fig. \ref{tad2}.


We now prove the following result:

\begin{theorem}
\label{th}
Tadfaces are not allowed by multi-orientability.
\end{theorem}
{\it Proof.} 
Let us first remark that, in order to prove this, one can forget about the ``middle'' strand since this strand never hooks to any of the ``external'' strands of an edge and thus cannot lead to a tadface.

\begin{figure}
\begin{center}
\includegraphics[scale=0.45,angle=0]{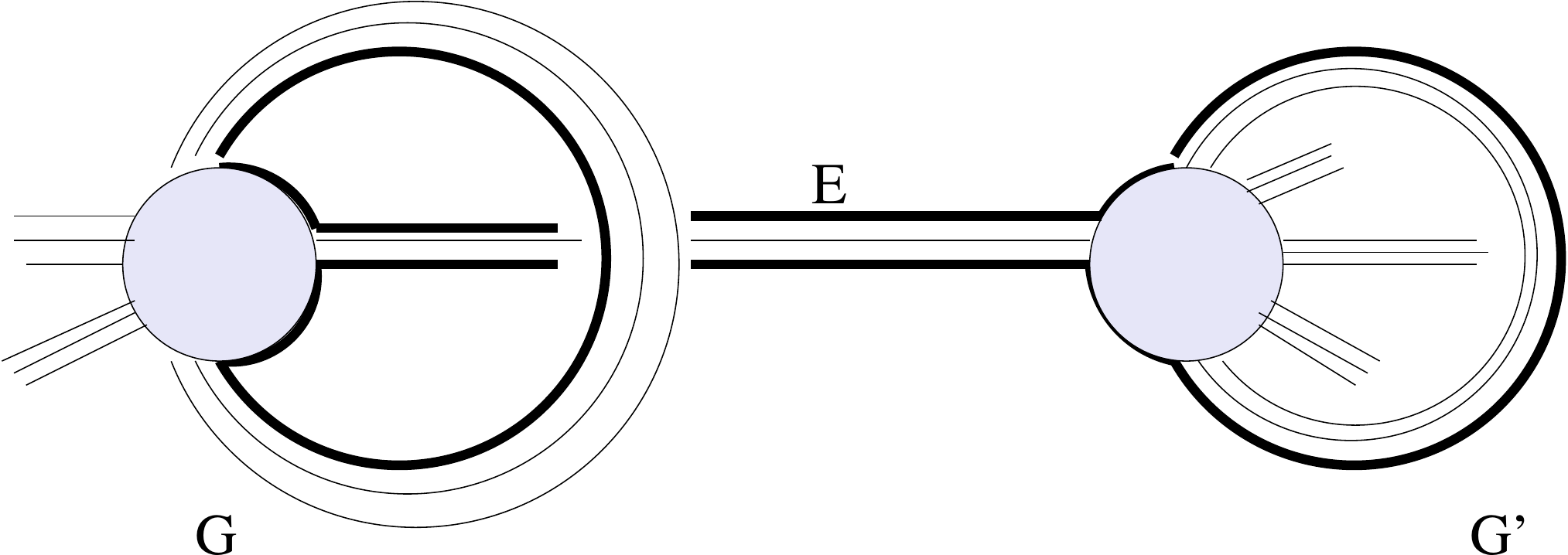}
\caption{A tadface in some general graph. One has an edge $E$ which is crossed twice by the respective circuit, once through one strand and once through another strand (nevertheless, the ``middle'' strand cannot be one of these two strands). The respective edge $E$ separates two subgraphs $G$ and $G'$ which has to be both irregular, the edge $E$ {\it alone} breaking a face for each of these subgraphs.}
\label{tf2}
\end{center}
\end{figure}

Furthermore, in order to obtain a tadface, one needs some  edge $E$ to be crossed twice when obtaining a face (according to Definition \ref{d:tf}) - see Fig. \ref{tf2}.
This means that one can identify two irregular subgraphs $G$ and $G'$ where the edge $E$ breaks a certain face of each of these subgraphs, while the other edges break other face(s) of the two subgraphs.

Joining together these two subgraphs through the edge $E$ leads to the tadface (see again Fig. \ref{tf2}).

We now prove the following intermediate result:

\begin{lemma}
\label{l}
One cannot have an irregular multi-orientable graph where one single edge $E$ breaks a face.
\end{lemma}
{\it Proof.} We suppose this graph exists.
Without any lose of generality we also suppose that the respective edge $E$ leaves from a ``-'' corner of some vertex. Then, one needs that a circuit leaves from a ``+'' corner and ends up in the opposite ``+'' corner of the same vertex (the respective circuit cannot close on a ``-'' corner because it would then leave inside another half-edge leaving from the missed ``+'' corner). The respective circuit can go through a certain number of vertices (even or odd); note however that the circuit hooks to a vertex through a ``-'' corner and needs to leave from a ``+'' corner without leaving half-edges inside. One can now directly check, using the parity constraints of the multi-orientable vertex, that such a circuit cannot be built, thus proving the Lemma. (QED).

\medskip

Before getting back to the proof of the main theorem, let us remark that irregular multi-orientable graphs (with four, six {\it etc.} external edges) do exist, but each broken face is broken by an even number of edges.
Using now  Lemma \ref{l} above, one concludes that the required subgraphs $G$ and $G'$ are forbidden by multi-orientability. This completes the proof. (QED)

\medskip

One can also obtain the following result, already announced in \cite{tadface}:

\begin{corollary}
Tadfaces are not allowed by colorability.
\end{corollary}
{\it Proof.} This follows directly from Theorem \ref{th} and Proposition \ref{inclu}.

\medskip

Let us now give {\it a posteriori} a physical motivation (leaving aside the pure combinatorial one) for getting interested to the presence or absence of these kind of tadfaces in GFT graphs. This motivation comes from the fact that it was proved in \cite{tf}  
that graphs without tadfaces  have better bounds per vertex within the framework of a BF model, like the ones we deal with in this paper.

\medskip

Furthermore, one  has:

\begin{corollary}
\label{c2}
``Non-planar'' generalized tadpoles are not allowed by multi-orientability.
\end{corollary}
{\it Proof.} A ``non-planar'' generalized tadpole is a graph with two external edges, each of these edges breaking a different face of the graph. One can now use Lemma \ref{l} to conclude that this type of graph is forbidden by multi-orientability. (QED)

\medskip 

As already announced above, this result was already known from the non-commutative QFT literature \cite{o1}.  

Finally, one has:

\begin{corollary}
\label{c3}
``Non-planar'' generalized tadpoles are not allowed by colorability.
\end{corollary}
{\it Proof.} This follows directly from Corollary \ref{c2} and Proposition \ref{inclu}. (QED)

\medskip

This last result on colorable models was already proved using a different method in the original paper \cite{color}.

\medskip

We resume all these results
 in the following table, which compares the colorable and the multi-orientable GFT models:

\begin{tabular}{ l |c| r }
 & colorable & multi-orientable \\
\hline\hline
  generalized planar tadpoles & forbidden & allowed\\ 
\hline
  generalized ``non-planar'' tadpoles & forbidden&forbidden \\
  
\hline
tadface & forbidden& forbidden
\end{tabular}

\bigskip

For the sake of completeness, let us also mention that the definition of bubbles or jackets of a general orientable graph, known from recent GFT literature (see for example \cite{fgo} for orientable graphs or \cite{matteo} for a general graph), is to be kept for the multi-orientable GFT introduced in this paper.

\medskip

Before ending this section, let us make one more comparison with the non-commutative QFT case. When dealing with models on the Moyal space (see again \cite{ncqft} and references within), orientability of the vertex discards the ``non-planar''-like tadpoles. As already mentioned in the Introduction, it is these two-point graphs which are mainly responsible for the appearance of the phenomenon of ultraviolet/infrared mixing. They have a $1/p^2$-like dependence in the external momenta (in the infrared regime of this one) and, when inserted into ``bigger'' graphs (which are thus non-planar), this external momenta becomes internal and needs to be integrated one leading to a new type of infrared divergence.

Nevertheless, some four-point graphs which are not discarded by the orientability of the vertex of the Moyal interaction still depend on the external momenta in 
a way which can lead to new types of divergences. This dependence is logarithmic and it leads indeed to the infrared divergence mentioned above when the respective four-point graph is similarly inserted into some ``bigger'' graph.

\begin{figure}
\begin{center}
\includegraphics[scale=0.5,angle=0]{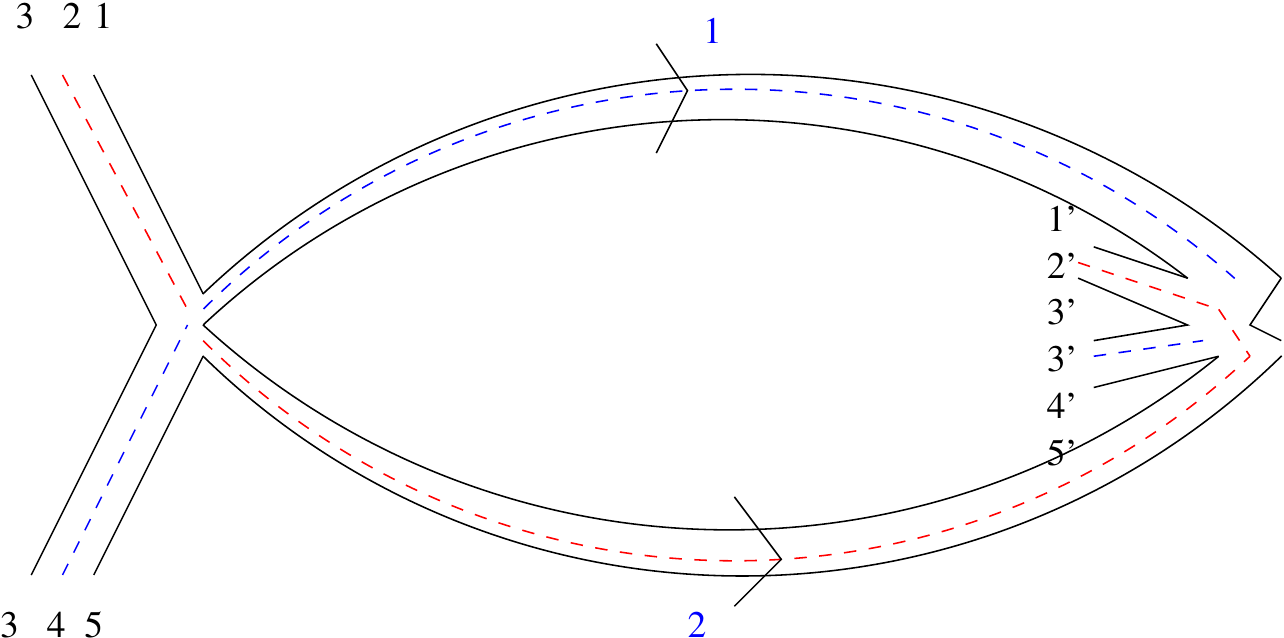}
\caption{An example of a non-colorable, multi-orientable four-point GFT graph. The indices in blue label the two internal edges. The rest of the indices refer to the group elements of the external edges.}
\label{4p-pli}
\end{center}
\end{figure}

When considering the GFT graphs, examples of such graphs are given in Fig. \ref{4p-pli}.
It is this type of graphs that we will carefully analyze in the next section.

\section{Feynman amplitude computations}
\renewcommand{\theequation}{\thesection.\arabic{equation}}
\setcounter{equation}{0}

We first investigate here the behavior of a non-colorable but multi-orientable graph, then that of a colorable (and multi-orientable) graph and finally, the one of a non-colorable and non-multi-orientable graph.


\subsection{A non-colorable, multi-orientable amplitude}

Let us now calculate the Feynman amplitude of the GFT graph  of Fig. \ref{4p-pli} We denote respectively by $h_{1}$ and by $h_{2}$ the two group elements associated to the internal edges $1$ and respectively $2$.

One has
\beqa
\label{1}
\int dh_1dh_2 \delta(g_1 h_1 h_2^{-1}g_5)\delta(g_4h_1g_{4'}^{-1})
\delta(g_2 h_2 g_{2'}^{-1})
\delta(g_{1'}^{-1}h_1^{-1}h_2g_{5'}^{-1}).
\eeqa
Performing the integral on $h_2$ using the third $\delta$ function in \eqref{1} and performing the integral on $h_1$ using the second $\delta$ function in \eqref{1} leads to the following result:
\beqa
\label{rez1}
\delta(g_1 g_{4'}g_4^{-1}g_{2'}^{-1}g_2g_5)
\delta(g_{1'}^{-1}g_4g_{4'}^{-1}g_2^{-1}g_{2'}g_{5'}^{-1}).
\eeqa
As expected the Feynman amplitude \eqref{1} is not divergent (this could have been directly stated from the fact there is no internal bubble of the GFT graph).

Nevertheless, an interesting phenomenon of some kind of ``ultraviolet/infrared'' mixing on the group manifold takes place. Thus, for
\beqa
\label{mix1}
g_1=g_5^{-1}g_2^{-1}g_{2'}g_4g_{4'}^{-1}g_1^{-1}
\mbox{ or } 
g_{1'}=g_4g_{4'}^{-1}g_2^{-1}g_{2'}g_{5'}^{-1}
\eeqa
the Feynman amplitude \eqref{rez1} becomes divergent. This comes from the fact that one has a non-trivial dependence of the amplitude on the external group momenta.

\subsection{A colorable, multi-orientable amplitude}

\begin{figure}
\begin{center}
\includegraphics[scale=0.5,angle=0]{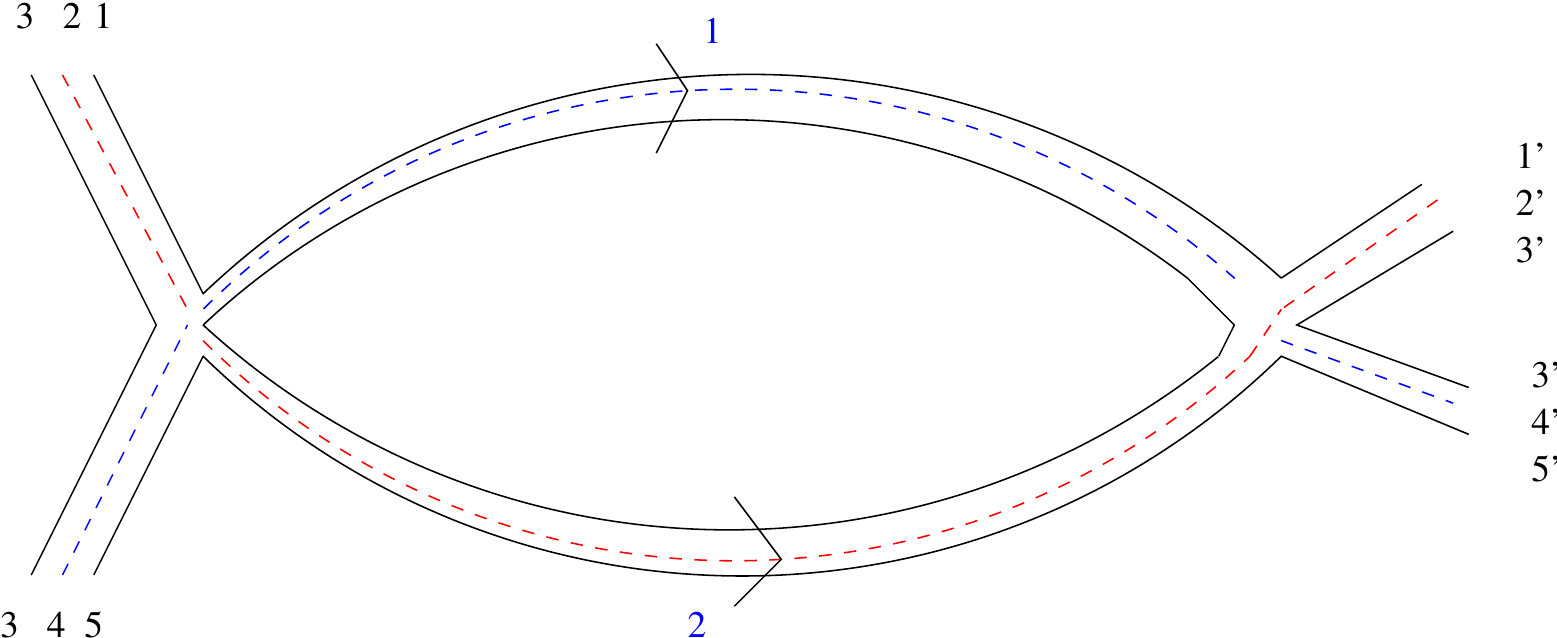}
\caption{An example of a colorable, multi-orientable four-point GFT graph. The indices in blue label the two internal edges. The rest of the indices refer to the group elements of the external edges.}
\label{4p-2}
\end{center}
\end{figure}

However, this phenomenon is not specific to the multi-orientable, non-colorable graphs graphs. In the case of the graph of Fig. \ref{4p-2}, 
which is colorable (and multi-orientable), the Feynman amplitude writes
\beqa
\label{2}
\int d h_1 dh_2\delta(g_1 h_1g_{1'}^{-1})\delta(g_4h_1g_{4'}^{-1})
\delta(h_1^{-1}h_2)
\delta(g_2 h_2 g_{2'}^{-1})\delta(g_5h_2g_{5'}^{-1}).
\eeqa
As above, we integrate on $h_2$ using the third $\delta$ function in \eqref{2} and we then integrate 
on $h_1$ using the first of the $\delta$ functions in \eqref{2}. The result writes
\beqa
\label{rez2}
\delta(g_4g_1^{-1}g_{1'}g_{4'}^{-1})
\delta(g_2g_1^{-1}g_{1'}g_{2'}^{-1})
\delta(g_5g_1^{-1}g_{1'}g_{5'}^{-1}).
\eeqa 
The group ``ultraviolet/infrared'' mixing described above is still present, several independent directions in the group
\beqa
\label{mix2}
g_{4'}=g_4g_1^{-1}g_{1'} 
\mbox{ or }
g_{2'}=g_2g_1^{-1}g_{1'}
\mbox{ or }
g_{5'}=g_5g_1^{-1}g_{1'},
\eeqa
turning the product \eqref{rez2} divergent.

\subsection{A non-colorable, non-multi-orientable amplitude}

\begin{figure}
\begin{center}
\includegraphics[scale=0.5,angle=0]{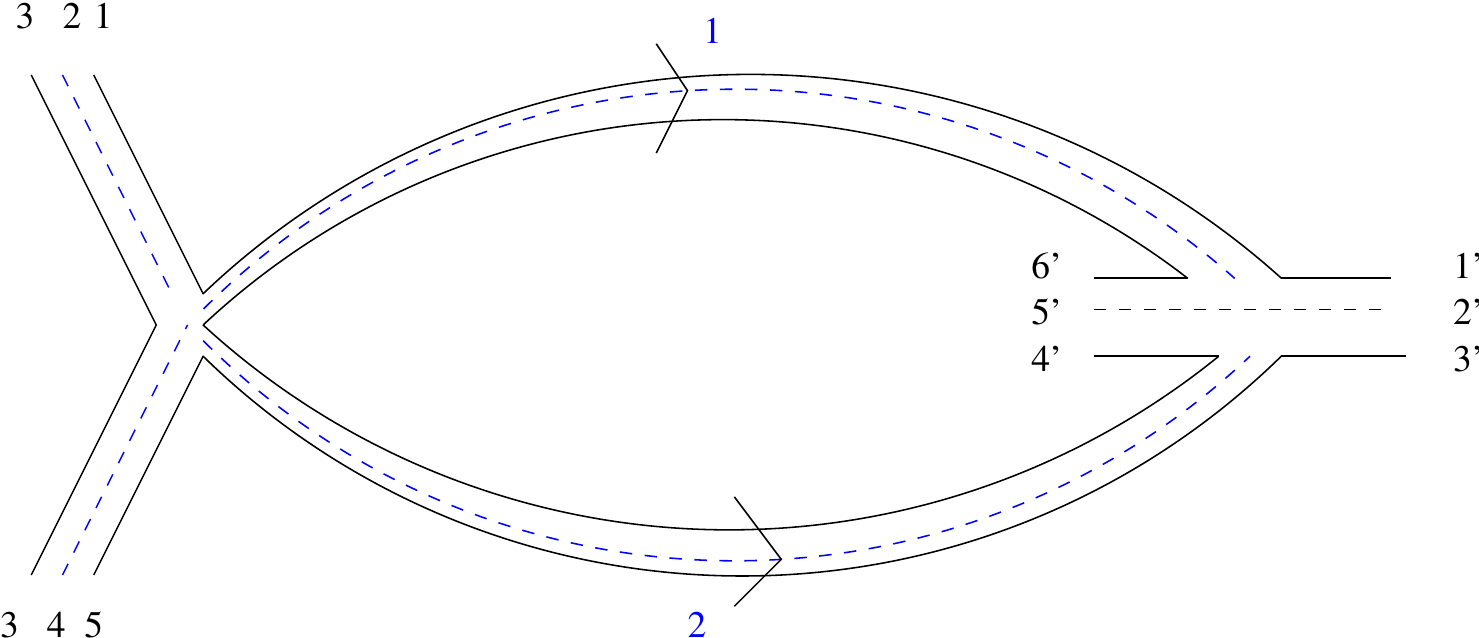}
\caption{An example of a non-colorable, non-multi-orientable four-point GFT graph. The indices in blue label the two internal edges. The rest of the indices refer to the group elements of the external edges.}
\label{4p-3}
\end{center}
\end{figure}

For the sake of completeness, we end this section by analyzing an associated non-colorable, non-multi-orientable GFT graph, like the one in Fig. \ref{4p-3}.
The Feynman amplitude of this graph writes
\beqa
\label{3}
\int d h_1 dh_2 \delta(g_1 h_1 g_{1'}^{-1})\delta(g_4h_1h_2^{-1}g_2)
\delta(g_{1'}h_1^{-1}h_2g_{4'})
\delta(g_2h_2h_1^{-1}g_4)\delta(g_5h_2g_{3'}^{-1}).
\eeqa
Integrating first on $h_1$ using the first of the $\delta$ functions in \eqref{3} and then on $h_2$ using the last of the $\delta$ function in \eqref{3} leads to the result:
\beqa
\delta (g_{6'}g_{1'}^{-1}g_1 g_5^{-1}g_{3'}g_{4'})
\delta(g_4g_1^{-1}g_{1'}g_{3'}^{-1}g_5 g_2)
\delta(g_2g_5^{-1}g_{3'}g_{1'}^{-1}g_1g_4).
\eeqa
This product of $\delta$ functions on external group elements 
can become infinite for any of the following 
independent group directions:
\beqa
g_{6'}=g_{4'}^{-1}g_{3'}^{-1}g_{5}g_1^{-1}g_{1'}
\mbox{ or }
g_4=g_2^{-1}g_5^{-1}g_{3'}g_{1'}^{-1}g_1
\mbox{ or }
g_2=g_4^{-1}g_1^{-1}g_{1'}g_{3'}^{-1}g_5.
\eeqa

The same type of phenomenon takes place when computing the Feynman amplitudes of the tadpoles of Fig. \ref{tad1} and respectively \ref{tad2}, which are non-colorable but multi-orientable and respectively non-colorable, non-multi-orientable.


\bigskip

Let us end this section with the following remark. The presence of the ``middle'' strand of these GFT edges (which makes the difference with respect to the combinatorial maps or ribbon maps of non-commutative QFT) is necessary for defining the bubbles, as it was already stated in \cite{io-JMP}. Moreover, this third strand is also required for the computation of Feynman amplitudes (which are related to the new concept of bubbles), as we have seen in detail in this section.

\section{Remark on multi-orientable GFT and renormalizability; comparison with the colorable GFT}
\renewcommand{\theequation}{\thesection.\arabic{equation}}
\setcounter{equation}{0}

In this section we give an example of a four-point GFT graph which, within the framework of the colorable models represents a quantum correction of a type not present in the bare action; within the framework of the multi-orientable models, this type of graph represents a quantum correction of a type already existent in the bare action. We then show a more general result regarding the quantum corrections of the two- and four-point functions within the multi-orientable framework.


\medskip

\begin{figure}
\begin{center}
\includegraphics[scale=0.3,angle=0]{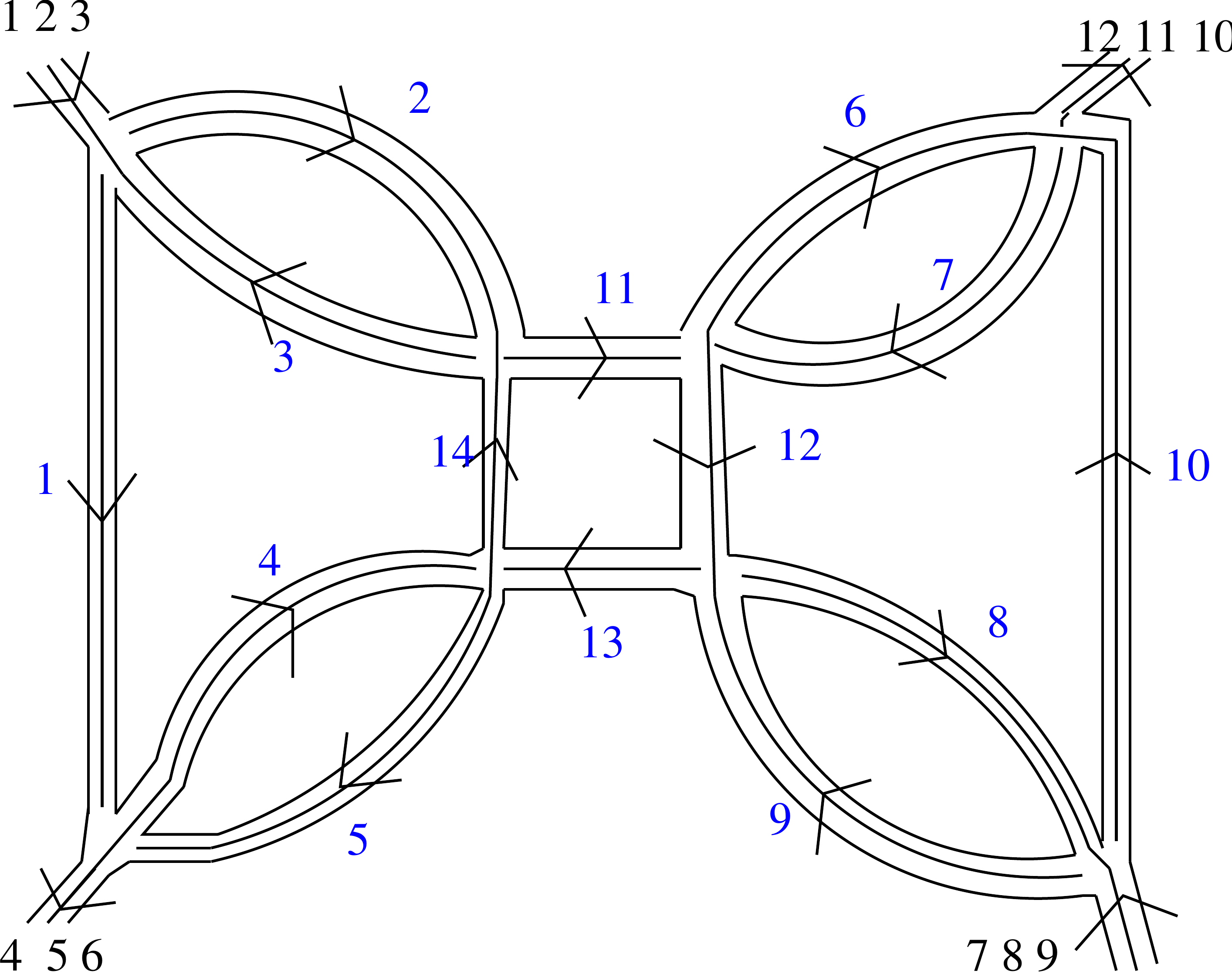}
\caption{An example of a four-point graph colorable (and multi-orientable) in three-dimensional GFT which, within the colorable framework, represents a quantum correction of a form not present in the bare action. In the multi-orientable framework, it represent a quantum correction of type $\bar \phi \phi\bar \phi \phi$, already present in the bare action. The indices in blue label the fourteen internal edges.}
\label{monstrul}
\end{center}
\end{figure}

Let us investigate the GFT behaviour of the tensor graph of Fig. \ref{monstrul}.
Denoting the group elements associated to the four external edges by $g_1,\ldots, g_{12}$ and the group elements associated to the fourteen internal edges by $h_1,\ldots,h_{14}$ as indicated in Fig. \ref{monstrul}, the Feynman amplitude writes:
\beqa
&&\int dh_1\ldots dh_{14}
\delta(g_1h_1g_4^{-1})\delta(h_1^{-1}h_2h_{14}^{-1}h_5)
\delta(h_1h_4h_{14}h_3)
\delta(g_3h_2h_{11}h_6g_{12}^{-1})\delta(h_2h_3)\nonumber\\
&&
\delta(g_2h_3^{-1}h_{11}h_7g_{11}^{-1})\delta(g_5^{-1}h_4h_{13}^{-1}h_8g_8)\delta(h_4h_5)
\delta(g_6h_5^{-1}h_{13}^{-1}h_9^{-1}g_7)
\delta(h_{14}h_{11}h_{12}h_{13})\nonumber\\
&&
\delta(h_6h_{10}^{-1}h_9h_{12}^{-1})\delta(h_6h_7)
\delta(h_7h_{12}h_8h_{10})
\delta(h_8h_9)
\delta(g_9h_{10}g_{10}^{-1}).
\eeqa
One remarks that we have a total nine internal $\delta$ functions ({\it i. e.} $\delta$ functions only in the internal parameters $h_1,\ldots, h_{14}$) which correspond to the nine internal faces of the tensor graph of Fig. \ref{monstrul}.
Upon direct inspection,
 one can then check that this graph is divergent as we did in the previous section (we recall that the Feynman amplitude computation is identical for Boulatov, colorable Boulatov or multi-orientable Boulatov models). Within the framework of the multi-orientable model \eqref{act:mo}, this represents a quantum correction of type 
$$ \bar \phi \phi \bar \phi \phi,$$
term which is present in the bare action. If one now colors the edges of this graph following the recipe indicated in section \ref{sec:3d}, the quantum correction here is for a term of type $\phi_p^4$ (where $p$ indices a specific color), which is not present in the bare colored action \eqref{act:col}.



\medskip

Let us now prove the following result:

\begin{theorem}
\noindent
\begin{itemize}
\item No two-point function of type $\phi\phi $ or $\bar\phi\bar\phi$ is permitted by multi-orientability.
\item No four-point function of a type distinct of $\phi\bar\phi\phi\bar\phi $
is permitted by multi-orientability.
\end{itemize}
\end{theorem}
{\it Proof.} Let us first prove the first point above. If one has two external edges then one has $B\le 2$. If $B=2$, these types of graphs are known to be not permitted by multi-orientability (see subsection $3.1$). If $B=1$, this means that one has the two external edges breaking the same (external) face. 
Using parity arguments like the ones developed for the proof of Theorem $3.1$, one can prove by induction (for example, first on the number of external vertices and then on the number of internal vertices) that such graphs cannot be constructed. This can be easily seen by direct inspection, since the number of edges to be added at an induction step is finite (namely, four).

Let us now prove the last item. As above, we suppose that a graph $G$ of type $\phi\phi\phi\bar\phi$ exists. Let us remark that this particular choice does not restrict the generality of the statement. This means that one can identify  a subgraph $G'$ of the original four-point graph $G$ such that two external edges of type $\phi$ and $\phi$ follows. This is however excluded and thus concludes the proof. (QED).

\medskip

Let us emphasize that the result above is valid at any order in perturbation theory.

\medskip

These phenomena appear as some kind of consequence of the 
fact that the 
restrictivity of the colorability condition is more significant then the restrictivity of the multi-orientability condition.

\section{Generalization to four dimensional GFT models}
\renewcommand{\theequation}{\thesection.\arabic{equation}}
\setcounter{equation}{0}

\begin{figure}
\begin{center}
\includegraphics[scale=0.25,angle=0]{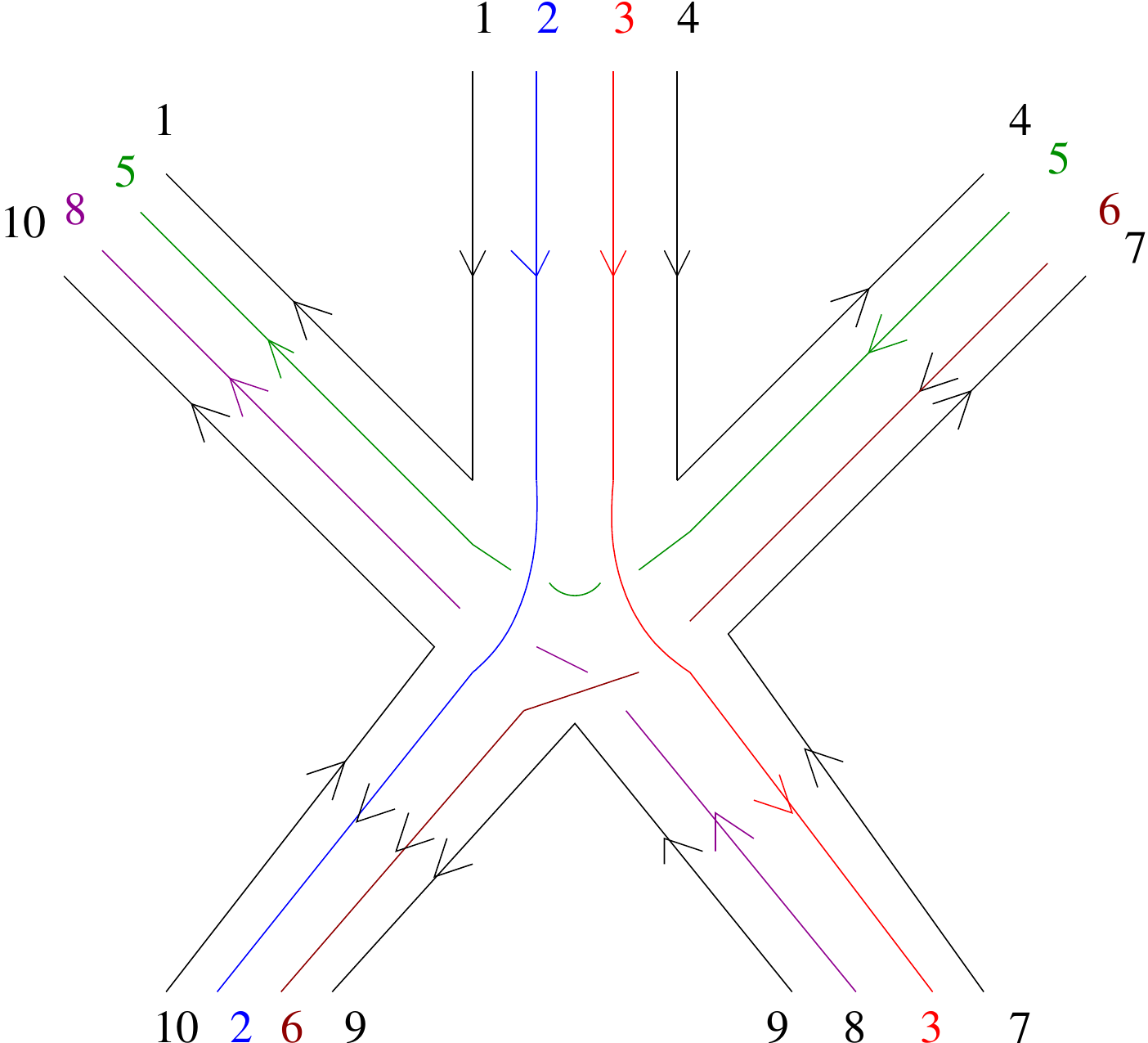}
\caption{Four-dimensional vertex for GFT models, topological or not.}
\label{vertex4d}
\end{center}
\end{figure}

\begin{figure}
\begin{center}
\includegraphics[scale=0.1, angle=0]{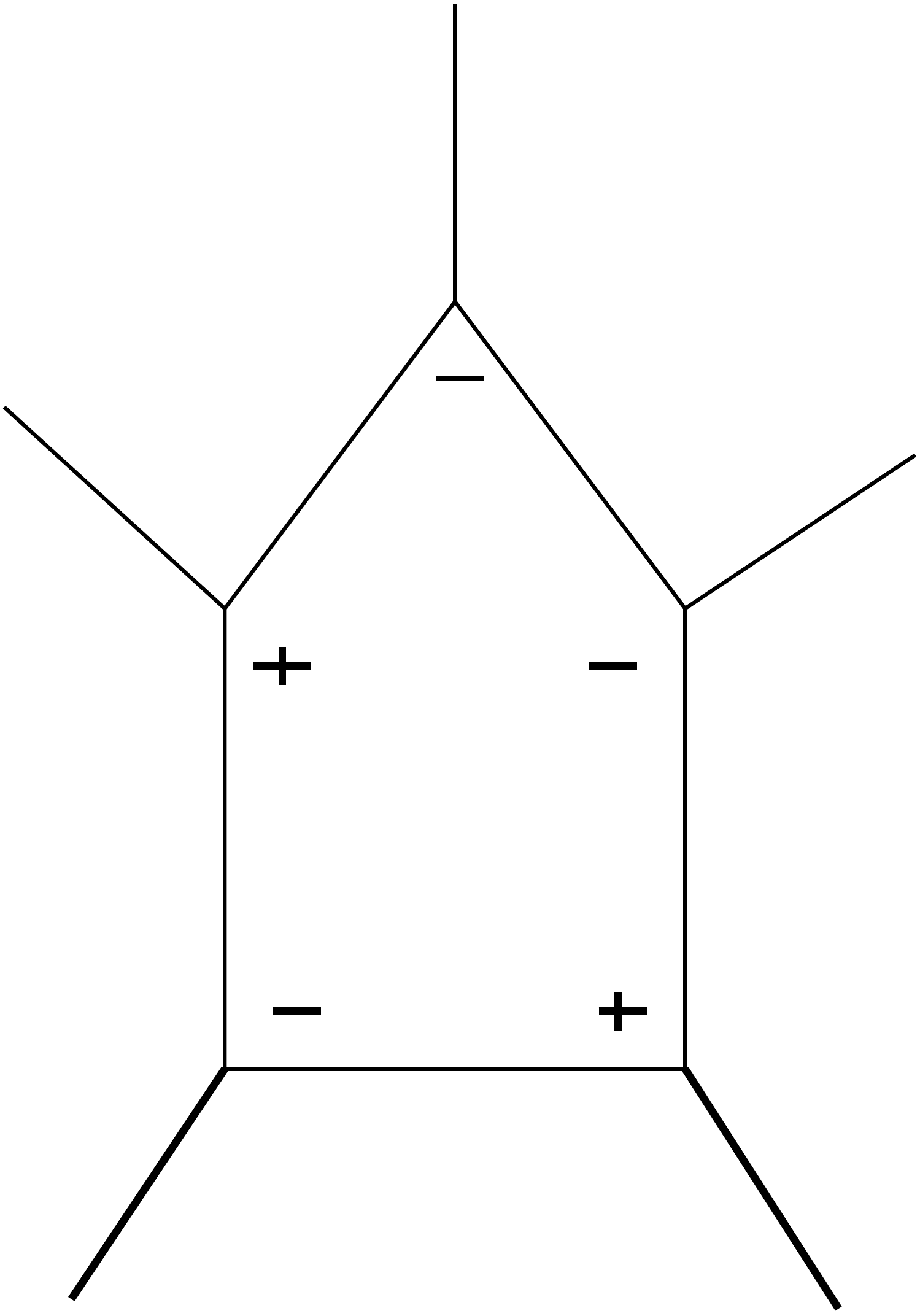}
\caption{A first possibility of an orientable vertex for four-dimensional GFT.}
\label{vertex1}
\end{center}
\end{figure}

\begin{figure}
\begin{center}
\includegraphics[scale=0.1,angle=0]{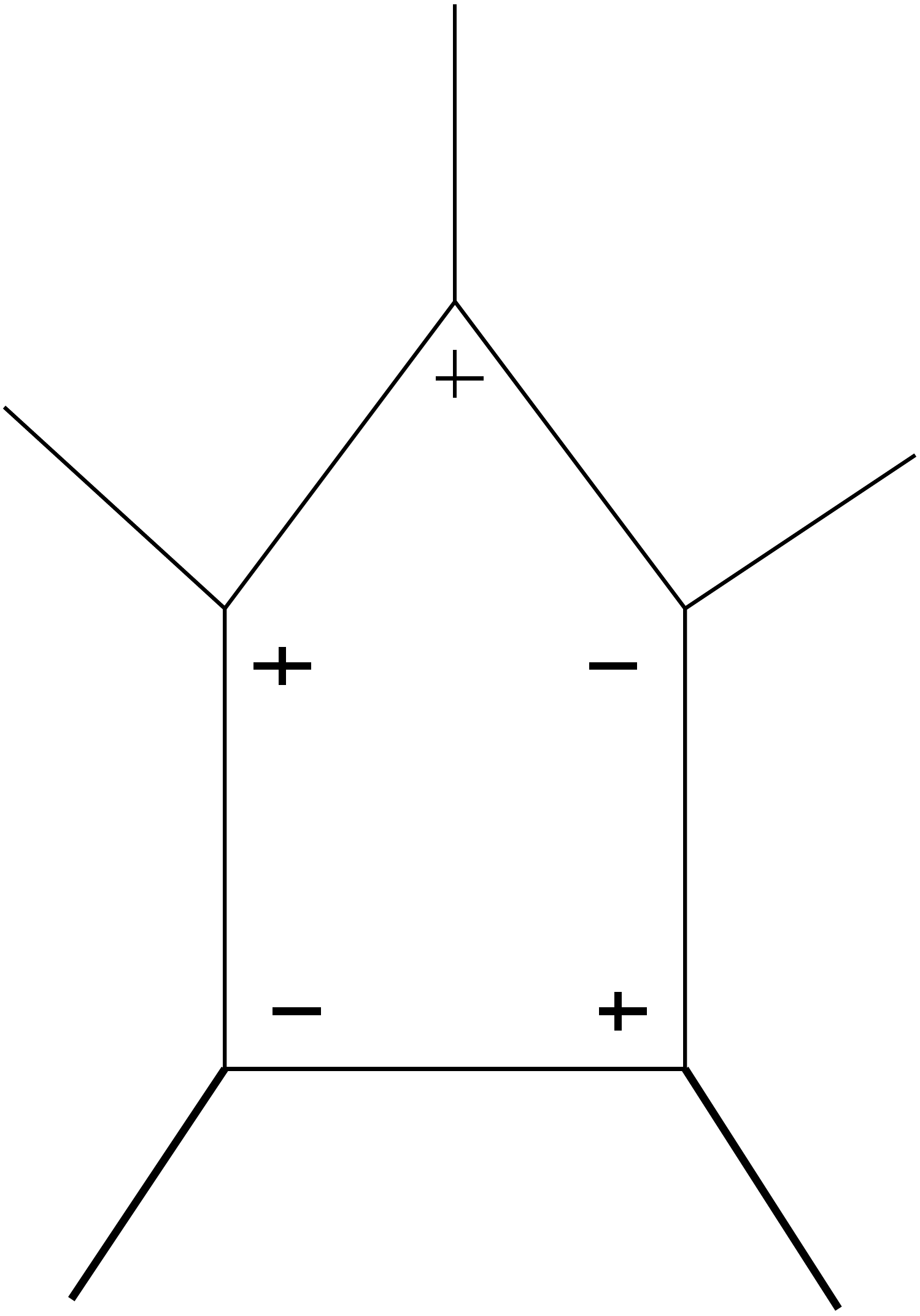}
\caption{A second possibility of an orientable vertex for four-dimensional GFT.}
\label{vertex2}
\end{center}
\end{figure}

The definition of the multi-orientable models to GFT in even dimension is not straightforward. This comes from the fact that the interaction $\phi^{D+1}$ is odd; one thus has two inequivalent choices of distribution of the $-$ and $+$ signs on the corners of the vertex. For the four-dimensional case (where the vertex is given in Fig. \ref{vertex4d}),
case of interest for quantum gravity, these two possibilities of interactions are given in Fig. \ref{vertex1} and \ref{vertex2}.

The action of the proposed model thus writes
\beqa
\label{act:mo4}
S[\phi]=\frac 12 \int \bar \phi \phi + \frac{\lambda_1}{5!}\int \bar \phi \phi \bar \phi \phi \bar \phi + \frac{\lambda_2}{5!}\int  \phi \bar \phi \phi \bar \phi \phi,
\eeqa
where, as in \eqref{act:col} or \eqref{act:mo}, the integrations over the group are left implicit. Moreover, 
in order to keep generality, we choose to have two distinct coupling constants $\lambda_1$ and $\lambda_2$; they are {\it a priori} free to flow under renormalization group.

The considerations of the previous sections extend for this four-dimensional multi-orientable GFT model.

\section{Conclusions and perspectives}
\renewcommand{\theequation}{\thesection.\arabic{equation}}
\setcounter{equation}{0}

We have  introduced  multi-orientable GFT in this paper, as a way of simplifying the topology and combinatorics of GFT; this simplification is different of the one proposed within the colorability of GFT. An analysis 
of the differences between 
the classes of tensor graphs discarded by these two types of models has been done and some Feynman amplitude computations have been performed. Moreover, we have given explicit examples of  two- and four-point graphs which, in the colorable framework lead to quantum corrections of types not present in the bare action, while in the multi-orientable framework lead to quantum correction of types already presented in the bare action. Finally, a generalization from three-dimensional to four-dimensional models has been proposed.

\medskip

Let us now give more details on the general usefulness of such a concept of multi-orientable GFT. As in the case of colored GFT, the first motivation comes from the fact that GFT are known to have an extremely involved combinatorial and topological structure \cite{not-ori2}. Multi-orientable GFT discards ({\it via} the indicated QFT recipe) some of these graphs and this means that explicit mathematical manipulation are more likely to be pushed further in this simpler case rather then in the general case. One can thus see multi-orientable (or colorable) GFT (at least) as a laboratory for testing various QFT tools before tackling the same problems in the general case.

Moreover, the results of the section $5$ show that multi-orientable GFT have an encouraging behaviour when one has in mind renormalizability studies. Thus, these models are shown to comport better then colorable models, where one has quantum corrections of types non present in the bare action.

Finally, let us also recall a general motivation for colored GFT (see again \cite{color}), motivation which also applies for multi-orientable GFT. This motivation comes from an analogy with matrix models, where only identically distributed models are topological in some scaling limit \cite{2d}. Nevertheless, it is topology which governs the power counting of more involved models (the renormalizable ones, \cite{gw} and \cite{gmrt}). It is this mathematical feature that the multi-orientable (and colorable) GFT is exploiting.

\medskip

Since this paper gives a proposal for a new type of GFT models, the perspectives for future work on this subject appears to us as particularly vast. Thus, it would be interesting to check whether or not the various achievements obtained within the framework of colorable
 GFT models (see again \cite{dev-cgft}, \cite{bgo}, \cite{ryan}, \cite{j2}, \cite{tf} and references within) can be adapted (and in what conditions) to  multi-orientable GFT. 

A first perspective may be the investigation of the issue of orientability of the piecewise linear pseudo-manifold associated to the graph
of the perturbative expansion of multi-orientable GFT. One needs to check weather or not the proof of the paper \cite{cara} for the colored Boulatov models adapts to the multi-orientable framework described here.

Among other possible perspectives, we can list here: the definition of some computable cellular homology, the (celebrated) $1/N$-expansion, the investigation of scaling behavior, of Borel summability or of (local and global) gauge transformations.

Regarding this last point of gauge symmetries, let us stress the following issue, which appears to be of crucial importance. In an ``ordinary'' gauge theory, for example in QED, one can redundantly integrate over a continuous infinity of physically equivalent configurations. This type of problem can be fixed using Fadeev-Popov-like tricks (see again textbooks like \cite{qft2}). 
This situation can also arise in GFT (colorable or not, multi-orientable or not) and it seems that one needs to make sure if this indeed happens or not in order to correctly interpret the divergences of the various models \cite{jacques}.

More precisely, in \cite{bgo} a global diffeomorphism symmetry was identified for the colored Boulatov model. The first step in the project mentioned above is thus to look for a local form of such a symmetry and then to check whether or not one can factorize (a part of) the colored Boulatov divergences. The same can be proposed for multi-orientable models. In this case, however, one has to look first for an equivalent of the diffeomorphism symmetry identified in \cite{bgo}.


\section*{Acknowledgments}
 The author acknowledges the CNRS PEPS grant ``CombGraph''. The grants PN 09 37 01 02 and CNCSIS Tinere Echipe 77/04.08.2010 are also acknowledged.


\begin{thebibliography}{99}


\bibitem{graph-books}
W. T. Tutte, ``Graph Theory'' , Reading, Mass.: Addison-Wesley, 1984. R. Diestel. ``Graph Theory''
2010. 
Springer-Verlag.

\bibitem{qft2}
C. Itzykson and J.-B. Zuber, "Quantum field theory", McGraw-Hill, 1980.
M. Peskins and D. Schroeder, "An Introduction to Quantum Field Theory", 1995, Addison-Wesley Advanced Book Program.

\bibitem{books-QFT}
V. Rivasseau, "From Perturbative to Constructive Renormalization", Princeton Series in Physics, 1991.
H. Kleinert, V. Schulte-Frohlinde, "Critical Properties of $\Phi^4$-Theories", World Scientific Publishing Co Pte Ltd.



\bibitem{io-SLC}
A.~Tanasa,
  ``Some combinatorial aspects of quantum field theory,''
    [arXiv:1102.4231 [math.CO]], submitted.

\bibitem{maps}
G. Schaeffert, Habilitation thesis, 2009, \'Ecole Polytechnique. G. Chapuy, ``Combinatoire bijective des cartes de genre sup\'erieur'', PhD thesis, \'Ecole Polytechnique (in French).

\bibitem{ncqft}
V.~Rivasseau,
  ``Non-commutative Renormalization,'' S\'eminaire Henri Poincar\'e,
    [arXiv:0705.0705 [hep-th]].
   A.~Tanasa,
  ``Translation-invariant noncommutative renormalization,''
  SIGMA {\bf 6}, 047 (2010).
  [arXiv:1003.4877 [hep-th]].
J.~Ben Geloun, A.~Tanasa,
``One-loop beta functions of a translation-invariant renormalizable noncommutative scalar model,''
  Lett.\ Math.\ Phys.\  {\bf 86}, 19-32 (2008).
  [arXiv:0806.3886 [math-ph]].
A.~Tanasa,
  ``Scalar and gauge translation-invariant noncommutative models,''
  Rom.\ J.\ Phys.\  {\bf 53}, 1207-1212 (2008).
  [arXiv:0808.3703 [hep-th]].
 A.~Tanasa, D.~Kreimer,
``Combinatorial Dyson-Schwinger equations in noncommutative field theory,''
J. Noncomm. Geom. (in press) 
  [arXiv:0907.2182 [hep-th]].


\bibitem{2d}
F. David, 
``A model of random surfaces with nontrivial critical behavior'', 
Nucl. Phys. {\bf B} 257 (1985). 

\bibitem{2d2}
P. Ginsparg, hep-th/9112013.
P. di Francesco,
 ``2D quantum gravity, matrix models and graph combinatorics,'' NATO Advanced Study Institute, 
   [math-ph/0406013].


\bibitem{gft}
  L.~Freidel,
  ``Group field theory: An overview,''
  Int.\ J.\ Theor.\ Phys.\  {\bf 44}, 1769 (2005)
  [arXiv:hep-th/0505016].
  D.~Oriti,
  ``The group field theory approach to quantum gravity,''
  arXiv:gr-qc/0607032.
  V.~Rivasseau,
  ``Towards Renormalizing Group Field Theory,''
  PoS {\bf CNCFG2010}, 004 (2010).
  [arXiv:1103.1900 [gr-qc]].

\bibitem{sf}
 M.~P.~Reisenberger, C.~Rovelli,
``Space-time as a Feynman diagram: The Connection formulation,''
  Class.\ Quant.\ Grav.\  {\bf 18}, 121-140 (2001).
  [gr-qc/0002095].
A.~Perez,
``Spin foam models for quantum gravity,''
  Class.\ Quant.\ Grav.\  {\bf 20}, R43 (2003)
  [gr-qc/0301113].

\bibitem{dia}
C.~Rovelli,
``A critical look at strings,''
  [arXiv:1108.0868 [hep-th]].
 C.~Rovelli,
``A Dialog on quantum gravity,''
  Int.\ J.\ Mod.\ Phys.\  {\bf D12}, 1509-1528 (2003).
  [hep-th/0310077].
D.~Oriti,
``Approaches to quantum gravity: Toward a new understanding of space, time and matter,''
  Cambridge, UK: Cambridge Univ. Pr. (2009) 583 p


\bibitem{gft-regain}
 J.~Magnen, K.~Noui, V.~Rivasseau, M. Smerlak,,
 ``Scaling behaviour of three-dimensional group field theory,''
  Class.\ Quant.\ Grav.\  {\bf 26 } (2009)  185012.
  [arXiv:0906.5477 [hep-th]].
  J.~Ben Geloun, V.~Bonzom,
``Radiative corrections in the Boulatov-Ooguri tensor model: The 2-point function,''
  Int.\ J.\ Theor.\ Phys.\  {\bf 50}, 2819-2841 (2011).
  [arXiv:1101.4294 [hep-th]].
    J.~Ben Geloun,
``Classical Group Field Theory,''
  [arXiv:1107.3122 [hep-th]].
  M.~Dupuis, F.~Girelli, E.~R.~Livine,
``Spinors and Voros star-product for Group Field Theory: First Contact,''
  [arXiv:1107.5693 [gr-qc]].
  V.~Bonzom, M.~Smerlak,
``Bubble divergences: sorting out topology from cell structure,''
  [arXiv:1103.3961 [gr-qc]].
A.~Baratin and D.~Oriti,
``Group field theory with non-commutative metric variables,''
  Phys.\ Rev.\ Lett.\  {\bf 105}, 221302 (2010)
  [arXiv:1002.4723 [hep-th]].
A.~Baratin and D.~Oriti,
``Ten questions on Group Field Theory (and their tentative answers),''
  arXiv:1112.3270 [gr-qc].

  \bibitem{uv/ir}
  S.~Minwalla, M.~Van Raamsdonk, N.~Seiberg,
  ``Noncommutative perturbative dynamics,''
  JHEP {\bf 0002}, 020 (2000).
  [hep-th/9912072].
J.~Magnen, V.~Rivasseau, A.~Tanasa,
  ``Commutative limit of a renormalizable noncommutative model,''
  Europhys.\ Lett.\  {\bf 86}, 11001 (2009).
  [arXiv:0807.4093 [hep-th]]

\bibitem{bahns}
D.~Bahns, S.~Doplicher, K.~Fredenhagen and G.~Piacitelli,
 ``Ultraviolet finite quantum field theory on quantum space-time,''
  Commun.\ Math.\ Phys.\  {\bf 237}, 221 (2003)
  [hep-th/0301100].

\bibitem{braided}
R.~Oeckl,
``Untwisting noncommutative R**d and the equivalence of quantum field theories,''
  Nucl.\ Phys.\ B {\bf 581}, 559 (2000)
  [hep-th/0003018].
Y.~Sasai and N.~Sasakura,
``Braided quantum field theories and their symmetries,''
  Prog.\ Theor.\ Phys.\  {\bf 118}, 785 (2007)
  [arXiv:0704.0822 [hep-th]].


\bibitem{color}
  R.~Gurau,
``Colored Group Field Theory,'' Commun.\ Math.\ Phys.\  {\bf 304}, 69-93 (2011).
  arXiv:0907.2582 [hep-th].

\bibitem{bipartite}
Weisstein, Eric W. "Bipartite Graph." From MathWorld--A Wolfram Web Resource. http://mathworld.wolfram.com/BipartiteGraph.html

\bibitem{cristal}
S. Lins, ``Gems, computers and attractors for 3-manifolds'', {\it Series on Knots and Everything}, {\bf 5}, World Scientific, 1995. M. Ferri and C. Gagliardi, ``Crystallisation moves'', {\it Pacific J. Math.}, {\bf 100}, 1, 1982. M. Pezzana, ``Sulla struttura topologica delle variet\`a compatte'', Atti Sem. Mat. Fis. Univ. Modena
23 (1974) 269-277.

\bibitem{cara}
 F.~Caravelli,
 ``A simple proof of orientability in the colored Boulatov model,''
  arXiv:1012.4087 [math-ph].


\bibitem{dev-cgft}
J.~Ben Geloun,
``Ward-Takahashi identities for the colored Boulatov model,''
  [arXiv:1106.1847 [hep-th]].
V.~Bonzom, R.~Gurau, A.~Riello, V.~Rivasseau,
``Critical behavior of colored tensor models in the large N limit,''
  [arXiv:1105.3122 [hep-th]].
R.~Gurau,
``The complete 1/N expansion of colored tensor models in arbitrary dimension,''
  [arXiv:1102.5759 [gr-qc]].
 R.~Gurau, V.~Rivasseau,
``The 1/N expansion of colored tensor models in arbitrary dimension,''
  [arXiv:1101.4182 [gr-qc]].

\bibitem{bgo}
 A.~Baratin, F.~Girelli, D.~Oriti,
  ``Diffeomorphisms in group field theories,''
  Phys.\ Rev.\  {\bf D83}, 104051 (2011).
  [arXiv:1101.0590 [hep-th]].


\bibitem{ryan}
J.~P.~Ryan,
``Tensor models and embedded Riemann surfaces,''
  [arXiv:1104.5471 [gr-qc]].

\bibitem{j2}
R.~Gurau,
  ``The 1/N expansion of colored tensor models,''
  Annales Henri Poincare {\bf 12}, 829-847 (2011).
  [arXiv:1011.2726 [gr-qc]].

\bibitem{tf}
J.~Ben~Geloun, J.~Magnen, V.~Rivasseau,
  ``Bosonic Colored Group Field Theory,'' Eur.\ Phys.\ J.\  {\bf C70}, 1119-1130 (2010).
   [arXiv:0911.1719 [hep-th]].





\bibitem{o1}
I.~Chepelev, R.~Roiban,
``Convergence theorem for noncommutative Feynman graphs and renormalization,''
  JHEP {\bf 0103}, 001 (2001).
  [hep-th/0008090].

\bibitem{orient-ncqft}
F.~Vignes-Tourneret,
``Renormalization of the Orientable Non-commutative Gross-Neveu Model,''
  Annales Henri Poincare {\bf 8}, 427-474 (2007).
  [math-ph/0606069].
 R.~Gurau, A.~Tanasa,
``Dimensional regularization and renormalization of non-commutative QFT,''
  Annales Henri Poincare {\bf 9}, 655-683 (2008).
  [arXiv:0706.1147 [math-ph]].
V.~Rivasseau, A.~Tanasa,
  ``Parametric representation of ''covariant'' noncommutative quantum field theory models,''
  Commun.\ Math.\ Phys.\  {\bf 279}, 355-379 (2008).
  [math-ph/0701034].
R.~Gurau, V.~Rivasseau,
``Parametric representation of noncommutative field theory,''
  Commun.\ Math.\ Phys.\  {\bf 272}, 811-835 (2007).
  [math-ph/0606030]
A.~Tanasa,
``Overview of the parametric representation of renormalizable non-commutative field theory,''
  J.\ Phys.\ Conf.\ Ser.\  {\bf 103}, 012012 (2008).
  [arXiv:0709.2270 [hep-th]].


\bibitem{boulatov}
D.~V.~Boulatov, Mod. Phys. Lett. {\bf  A7} (1992) 1629;
eprint {\tt hep-th/9205090}.

\bibitem{ooguri}
H.~Ooguri,
``Topological lattice models in four-dimensions,''
 Mod.\ Phys.\ Lett.\  {\bf A7}, 2799-2810 (1992).
 [hep-th/9205090].



\bibitem{io-5}
T.~Krajewski, J.~Magnen, V.~Rivasseau, A.~Tanasa, P.~Vitale,
``Quantum Corrections in the Group Field Theory Formulation of the EPRL/FK Models,''
  Phys.\ Rev.\  {\bf D82}, 124069 (2010).
  [arXiv:1007.3150 [gr-qc]].

\bibitem{not-ori}
  R.~De Pietri, L.~Freidel, K.~Krasnov, C.~Rovelli,
``Barrett-Crane model from a Boulatov-Ooguri field theory over a homogeneous space,''
  Nucl.\ Phys.\  {\bf B574}, 785-806 (2000).
  [hep-th/9907154].

\bibitem{not-ori2}
R.~De Pietri, C.~Petronio,
``Feynman diagrams of generalized matrix models and the associated manifolds in dimension 4,''
  J.\ Math.\ Phys.\  {\bf 41}, 6671-6688 (2000).
  [gr-qc/0004045].


\bibitem{fusy}
Eric Fusy, private communication.

\bibitem{tadface}
J.~Ben Geloun, T.~Krajewski, J.~Magnen, V.~Rivasseau,
``Linearized Group Field Theory and Power Counting Theorems,''
  Class.\ Quant.\ Grav.\  {\bf 27}, 155012 (2010).
  [arXiv:1002.3592 [hep-th]].

\bibitem{fgo}
  L.~Freidel, R.~Gurau and D.~Oriti,
  ``Group field theory renormalization - the 3d case: power counting of
  divergences,''
  Phys.\ Rev.\  D {\bf 80}, 044007 (2009)
  [arXiv:0905.3772 [hep-th]].

\bibitem{matteo}
V.~Bonzom, M.~Smerlak,
``Bubble divergences from cellular cohomology,''
  Lett.\ Math.\ Phys.\  {\bf 93}, 295-305 (2010).
  [arXiv:1004.5196 [gr-qc]].

\bibitem{io-JMP}
A.~Tanasa,
  ``Generalization of the Bollob\'as-Riordan polynomial for tensor graphs,'' J. Math. Phys. (in press)
    [arXiv:1012.1798 [math.CO]].






\bibitem{gw}
 H.~Grosse and R.~Wulkenhaar,
 ``Renormalization of phi**4 theory on noncommutative R**4 to all orders,''
  Lett.\ Math.\ Phys.\  {\bf 71}, 13 (2005)
  [hep-th/0403232].

\bibitem{gmrt}
 R.~Gurau, J.~Magnen, V.~Rivasseau and A.~Tanasa,
``A Translation-invariant renormalizable non-commutative scalar model,''
  Commun.\ Math.\ Phys.\  {\bf 287}, 275 (2009)
  [arXiv:0802.0791 [math-ph]].

\bibitem{jacques}
J. Magnen, private communication.

\end{thebibliography}
\end{document}